\RequirePackage{fix-cm}
\documentclass[smallextended]{svjour3}       
\smartqed  
\usepackage[utf8]{inputenc}
\usepackage{amsmath}
\usepackage{amsfonts}
\usepackage{amssymb}
\usepackage{enumerate}
\usepackage{mathtools}

\providecommand{\norm}[1]{\lVert#1\rVert}
\journalname{Calculus of Variations and PDEs}

\begin{document}

\title{Optimal transport of vector measures\thanks{The author wishes to thank Bo'az Klartag for proposing to work on this problem and for useful discussions. The financial support of St John's College in Oxford, Clarendon Fund and EPSRC is gratefully acknowledged.  Part of this research was completed in Fall 2017 while the author was member of the Geometric Functional Analysis and Application program at MSRI, supported by the National Science Foundation under Grant No. 1440140. This research was also partly supported by the ERC Starting Grant 802689 CURVATURE.}
}

\author{Krzysztof J. Ciosmak}

\institute{Krzysztof J. Ciosmak\at
              Mathematical Institute and St John's College, University of Oxford, Oxford, United Kingdom\\
              \email{ciosmak@maths.ox.ac.uk, krzysztof.ciosmak@sjc.ox.ac.uk}} 
\date{Received: date / Accepted: date}

\maketitle

\begin{abstract}
We develop and study a theory of optimal transport for vector measures. We resolve in the negative a conjecture of Klartag, that given a vector measure on Euclidean space with total mass zero, the mass of any transport set is again zero. We provide a counterexample to the conjecture. We generalise the Kantorovich--Rubinstein duality to the vector measures setting. Employing the generalisation, we answer the conjecture in the affirmative provided there exists an optimal transport with absolutely continuous marginals of its total variation.


\end{abstract}

\section{Introduction}

In this note we develop theory of optimal transport of vector measures. Let us first briefly describe the topic of classical optimal transport.

\subsection{Optimal transport}\label{ssec:ot}

In 1781 Gaspard Monge (see \cite{Monge}) asked the following question: given two probability distributions $\mu,\nu$ on a metric space $(X,d)$, how to transfer one distribution onto the other in an optimal way. The criterion of optimality was to minimise the average transported distance. Since then the topic has been developed extensively and much of this development has been done recently. We refer the reader to the books of Villani (see \cite{Villani2} and \cite{Villani1}) and to the lecture notes of Ambrosio (see \cite{Ambrosio3}) for a thorough discussion, history and applications of the optimal transport problem. 

The modern mathematical treatment of the problem has been initiated in 1942 by Kantorovich \cite{KantorovichURSS}, \cite{Kantorovich}. He proposed to consider a relaxed problem of optimising
\begin{equation*}
\int_{X\times X}d(x,y)d\pi(x,y)
\end{equation*}
among all transference plans $\pi$ between $\mu$ and $\nu$, i.e., the set $\mathrm{\Pi}(\mu,\nu)$ of Borel probability measures on $X\times X$ with respective marginal distributions equal to $\mu$ and to $\nu$.
The existence of an optimal transference plan is a straightforward consequence of the Prokhorov's theorem, provided that $X$ is separable.

The main question that has attracted a lot of attention is whether there exists an optimal transport plan, i.e., a Borel map $T\colon X\to X$ such that $T_{\#}\mu=\nu$ and the integral
\begin{equation*}
\int_Xd(x,T(x))d\mu(x)
\end{equation*}
is minimal.
If we knew that an optimal transference plan is concentrated on a graph of a Borel measurable function then we could infer the existence of an optimal transport plan. The first complete answer on Euclidean space, under regularity assumptions on the considered measures, was presented in a seminal paper \cite{Evans-Gangbo} of Evans and Gangbo. However, before that, Sudakov in \cite{Sudakov} presented a solution of the problem that contained a flaw. The flaw has been remedied by Ambrosio in \cite{Ambrosio3} and later by Trudinger and Wang in \cite{Trudinger} for the Euclidean distance and by Caffarelli, Feldman and McCann in \cite{Caffarelli} for distances induced by norms that satisfy certain smoothness and convexity assumptions. In \cite{Caravenna1} Caravenna has carried out the original strategy of Sudakov for general strictly convex norms and eventually Bianchini and Daneri in \cite{Bianchini1} accomplished the plan of a proof of Sudakov for general norms on finite-dimensional normed spaces. 

Let us describe briefly the strategy of Sudakov in the context of Euclidean spaces. We assume that the two Borel probability measures $\mu,\nu$ on $\mathbb{R}^n$ are absolutely continuous with respect to the Lebesgue measure. 

Let us recall that the paramount Kantorovich--Rubinstein duality formula tells that
\begin{equation}\label{eqn:kr}
\sup\Big\{\int_{\mathbb{R}^n}u d(\mu-\nu)\mid u\text{ is }1\text{-Lipschitz}\Big\}
\end{equation}
is equal to 
\begin{equation}\label{eqn:krrhs}
\inf\Big\{\int_{\mathbb{R}^n\times\mathbb{R}^n}\norm{x-y}d\pi(x,y)\mid \pi\in\mathrm{\Pi}(\mu,\nu)\Big\}.
\end{equation}
Let us take an optimal $u$ and an optimal $\pi$ in the two above optimisation problems. We may infer that 
\begin{equation}\label{eqn:mb}
u(x)-u(y)=\norm{x-y}\text{ for }\pi\text{-almost every }(x,y)\in X\times X.
\end{equation}
Consider the maximal sets on which $u$ is an isometry, called the \emph{transport rays}. We see that all transport has to occur on these sets. 
 Careful analysis of the Lipschitz function $u$ shows that the transport rays form a foliation of the underlying space $\mathbb{R}^n$ into line segments, up to Lebesgue measure zero. Moreover, the so-called mass balance condition holds true. This is to say, for any Borel set $A$ that is a union of some collection of transport rays there is $\mu(A)=\nu(A)$; see e.g. \cite{Evans-Gangbo}.
Using the mass balance condition, we may construct an optimal transport by gluing together optimal maps for each of the transport rays; see e.g. \cite{Ambrosio3}.

This is one of the important observations that is employed in the localisation technique, which allows to reduce the dimension of a considered problem; see a paper of Klartag \cite{Klartag} for application of the technique to weighted Riemannian manifolds satisfying the curvature-dimension condition in the sense of Bakry and \'Emery \cite{Bakry1}, \cite{Bakry} and papers of Cavalletti, Mondino \cite{Cavalletti3}, \cite{Cavalletti2} for application in the setting of metric measure spaces. The localisation technique stems from convex geometry, but its generalisations have been employed to prove many novel results concerning functional inequalities, e.g. isoperimetric inequality in the metric measure spaces satisfying the synthetic curvature-dimension condition (see \cite{Cavalletti3}, \cite{Cavalletti2}). The latter notion was introduced in the foundational papers by Sturm \cite{Sturm1}, \cite{Sturm2} and by Lott and Villani \cite{Villani3} and allowed for development of a far-reaching, vast theory of metric measure spaces. We refer the reader to \cite{Ciosmak2} and references therein for a broader description of the localisation technique and its history.

\subsection{Optimal transport of vector measures}

The purpose of this article is to investigate multi-dimensional generalisation of the optimal transport problem and its connections with the localisation technique, as proposed by Klartag in \cite[Chapter 6]{Klartag}.

We  shall consider finite-dimensional linear spaces equipped with Euclidean norm and $1$-Lipschitz maps $u\colon\mathbb{R}^n\to\mathbb{R}^m$. 
A leaf $\mathcal{S}$ of a $1$-Lipschitz map $u\colon\mathbb{R}^n\to\mathbb{R}^m$ is a maximal set, with respect to the order induced by inclusion, such that the restriction $u|_{\mathcal{S}}$ is an isometry. This is to say, $\mathcal{S}$ is a leaf, whenever for any $x,y\in\mathcal{S}$ there is 
\begin{equation*}
\norm{u(x)-u(y)}=\norm{x-y}
\end{equation*}
and for any $z\notin\mathcal{S}$ there exists $x\in\mathcal{S}$ such that
\begin{equation*}
\norm{u(x)-u(z)}<\norm{x-z}.
\end{equation*}
The notion of leaves is a multi-dimensional generalisation of the notion of transport rays, see Subsection \ref{ssec:ot}, of the one-dimensional optimal transport theory. We refer the reader to \cite{Ciosmak2} for a thorough study of such leaves. Let us mention that such leaves form a convex partition of $\mathbb{R}^n$, up to Lebesgue measure zero. Moreover, any two such leaves may intersect only by their relative boundaries; see \cite{Ciosmak2} for the proofs.

Suppose now that we are given a Borel probability measure $\mu$ on $\mathbb{R}^n$, absolutely continuous with respect to the Lebesgue measure, that satisfies $m$ linear constrains. This is to say,
\begin{equation}\label{eqn:const}
\int_{\mathbb{R}^n}fd\mu=0
\end{equation}
for some integrable function $f\colon\mathbb{R}^n\to\mathbb{R}^m$ with finite first moments, i.e.,
\begin{equation*}
\int_{\mathbb{R}^n}\norm{f(x)}\norm{x}d\mu(x)<\infty.
\end{equation*}
Let $u\colon\mathbb{R}^n\to\mathbb{R}^m$ be a $1$-Lipschitz map such that
\begin{equation}\label{eqn:maxx}
\int_{\mathbb{R}^n}\langle u,f\rangle d\mu=\sup\Big\{\int_{\mathbb{R}^n}\langle v,f\rangle d\mu\mid v\colon\mathbb{R}^n\to\mathbb{R}^m \text{ is }1\text{-Lipschitz}\Big\}.
\end{equation}
Existence of $u$ follows by the Arzel\`a--Ascoli theorem.

A Borel subset $A$ of $\mathbb{R}^n$ shall be called a \emph{transport set} associated to $u$, whenever for any $x\in A$ that belongs to a unique leaf of $u$ and any $y\in\mathbb{R}^n$ such that 
\begin{equation*}
\norm{x-y}=\norm{u(x)-u(y)},
\end{equation*}
there is $y\in A$. In other words, a transport set is a Borel union of a collection of leaves of $u$.
In \cite[Chapter 6]{Klartag} it is conjectured that  for any transport set $A$ of $u$
\begin{equation}\label{eqn:mass}
\int_A fd\mu=0.
\end{equation}
This is a generalisation of the mass balance condition, mentioned in Subsection \ref{ssec:ot}.  The affirmative answer to the conjecture would imply that one may decompose any Borel probability measure $\mu$, satisfying $m$ linear constraints of the form (\ref{eqn:const}), into a mixture of measures, concentrated on pairwise disjoint convex subsets of $\mathbb{R}^n$ of dimension at most $m$, satisfying the same linear constraints; see also \cite{Ciosmak2} for a discussion of the decomposition.

If $m=1$ then (\ref{eqn:maxx}) is precisely the dual problem to the optimal transport problem for measures  $\rho_1, \rho_2$ given by formulae $d\rho_1=f_+d\mu$ and $d\rho_2=f_-d\mu$. As we see in (\ref{eqn:kr}), the dual problem, depends merely on the difference of measures, and therefore, it makes sense to consider the optimal transport for signed measures with total mass zero. 

Inspired by this observation, in Section \ref{sec:transport} we develop a theory of optimal transport with metric cost of vector measures of total mass zero and study its basic properties. The r\^ole of a vector measure in the problem considered above is played by the measure with density $f$ with respect to the measure $\mu$.

The precise formulation of the optimal transport problem for an $\mathbb{R}^m$-valued measure $\eta$ on a metric space $(X,d)$ that we deal with is as follows:
\begin{equation}\label{eqn:vectortr}
\inf\Big\{\int_{X\times X} d(x,y) d\norm{\pi}(x,y)\mid \mathrm{P}_1\pi-\mathrm{P}_2\pi=\eta\Big\}.
\end{equation}
Here $\mathrm{P}_1\pi$ and $\mathrm{P}_2\pi$ stand for the first and the second marginal of the $\mathbb{R}^m$-valued measure $\pi$ respectively. The assumption on $\eta$ is that 
\begin{equation*}
\int_Xd(x,x_0)d\norm{\mu}(x)<\infty\text{ for some }x_0\in X\text{ and }\eta(X)=0.
\end{equation*}
The above problem for $m=1$ simplifies to the original optimal transport problem, as follows readily by the Kantorovich--Rubinstein formula. We prove that for $m>1$ an analogue of this formula holds with (\ref{eqn:kr}) replaced by 
\begin{equation}\label{eqn:lip}
\sup\Big\{\int_{X}\langle u,d\eta\rangle\mid u\colon X\to\mathbb{R}^m\text{ is }1\text{-Lipschitz}\Big\}
\end{equation}
and with (\ref{eqn:krrhs}) replaced by (\ref{eqn:vectortr}). This is a content of Theorem \ref{thm:lip}. We also develop a theory of the Wasserstein space $\mathcal{W}(X,\mathbb{R}^m)$ of vector-valued measures. We identify its dual space as the space of vector-valued Lipschitz maps; see Theorem \ref{thm:reciprocals}. Theorem \ref{thm:ae} provides an analogue of (\ref{eqn:mb}) in the new setting. 

The conjecture of Klartag (see \cite[Chapter 6]{Klartag}) in the language of our theory of optimal transport of vector measures may be restated as follows. Suppose that we are given a vector measure $\mu$ on $\mathbb{R}^n$, with $\mu(\mathbb{R}^n)=0$, which is absolutely continuous with respect to Lebesgue measure. Let $u\colon\mathbb{R}^n\to\mathbb{R}^m$ be a $1$-Lipschitz map, with respect to Euclidean norms, that attains the supremum
\begin{equation}\label{eqn:lipr}
\sup\Big\{\int_{\mathbb{R}^n}\langle v,d\mu\rangle\mid v\colon\mathbb{R}^n\to\mathbb{R}^m\text{ is }1\text{-Lipschitz}\Big\}.
\end{equation} 
It is claimed in \cite{Klartag} that the following mass balance condition holds true
\begin{equation}\label{eqn:massbalance}
\mu(A)=0\text{ for any Borel set }A\text{ that is a union of a family of leaves of }u.
\end{equation}
Using the developed theory, in Section \ref{sec:counter} we resolve the conjecture in the affirmative, provided that there exists an optimal transport with  marginals of its total variation that are absolutely continuous with respect to the Lebesgue measure; see Theorem \ref{thm:balance}. Note that in the one-dimensional setting, the existence of such optimal transport is clear; see (\ref{eqn:kr}) and (\ref{eqn:krrhs}).

We provide a counterexample to the conjecture, for the case $m>1$; see Theorem \ref{thm:nonzero}. It shows that, in general, the mass balance condition (\ref{eqn:massbalance}) fails to be true. It follows that it may happen that an optimal transport with absolutely continuous marginals do not exist, unlike in the one-dimensional case. 

More generally, let $\mathcal{F}$ be any subset of $1$-Lipschitz maps that is locally uniformly closed. We prove that the mass balance condition (\ref{eqn:massbalance}) fails to be true, even when the variational problem (\ref{eqn:lipr}) is replaced by 
\begin{equation}
\sup\Big\{\int_{\mathbb{R}^n}\langle v,d\mu\rangle \mid v\in\mathcal{F}\Big\},
\end{equation}
unless $\mathcal{F}$ is trivial, i.e. consists merely of affine maps. This is shown for also any norm on $\mathbb{R}^n$ and any strictly convex norm on $\mathbb{R}^m$; see Theorem \ref{thm:affine}. 

Note that the outline of a proof of the conjecture suggested in \cite{Klartag} has a gap, as follows by the results of \cite{Ciosmak}.

Let us mention here that in \cite{Ciosmak2} the generalisation of the localisation technique to multiple constraints is studied. In there, a partition associated to any $1$-Lipschitz map $u\colon\mathbb{R}^n\to\mathbb{R}^m$, $m\leq n$, is studied thoroughly. It is established that any log-concave measure on $\mathbb{R}^n$ may be disintegrated with respect to this partition and that the resulting conditional measures, associated to leaves of maximal dimension, are again log-concave.
This result is also presented in the context of spaces satisfying the curvature-dimension condition $CD(\kappa,N)$, thus partially confirming another conjecture of Klartag \cite[Chapter 6]{Klartag}.

Let us also mention the existence of another approach to optimal transport of vector measures that differs from ours developed by Chen, Georgiou, Tannenbaum, Tyu, Li, Osher, Haber, Yamamoto (see \cite{Chen1}, \cite{Chen3} and \cite{Chen2}).

\subsection{Outline of the paper}

Section \ref{sec:transport} is devoted to development and study of the optimal transport theory of vector measures. We define a Wasserstein space and in Theorem \ref{thm:reciprocals} we identify its dual. Theorem \ref{thm:lip} provides an analogue of the Kantorovich--Rubinstein duality formula.

In Section \ref{sec:counter} we study the mass balance condition for vector measures. In Theorem \ref{thm:balance} we answer in the affirmative the a conjecture of Klartag, provided there exists an optimal transport with absolutely continuous marginals of its total variation. In Theorem \ref{thm:nonzero} we provide a counterexample to the conjecture, in the Euclidean setting.
In Theorem \ref{thm:affine} we resolve the conjecture in the negative in the general setting.

\section{Optimal transport of vector measures}\label{sec:transport}

In this section we develop the theory of optimal transport of vector measures.

Let $X$ be a metric space with metric $d$. Let $\mu$ be $\mathbb{R}^m$-valued Borel measure on $X$. If $\pi$ is a $\mathbb{R}^m$-valued Borel measure on $X\times X$, we write $\mathrm{P}_1\pi$ for the first \emph{marginal} of $\pi$, i.e. the measure given by 
\begin{equation*}
\mathrm{P}_1\pi(A)=\pi(A\times X),
\end{equation*}
for all Borel $A\subset X$, and $\mathrm{P}_2\pi$ for the second \emph{marginal} of $\pi$, 
\begin{equation*}
\mathrm{P}_2\pi(B)=\pi(X\times B),
\end{equation*}
for all Borel $B\subset X$. We shall consider a variational problem
\begin{equation}\label{eqn:KR}
\mathcal{I}(\mu)=\inf\bigg\{{\int_{X\times X}}d(x,y) d\norm{\pi}(x,y)\mid \pi\in \Gamma(\mu) \bigg\}.
\end{equation}
Here $\Gamma(\mu)$ is the set of all $\mathbb{R}^m$-valued Borel measures $\pi$ on $X\times X$ such that 
\begin{equation*}
\mu=\mathrm{P}_1\pi -\mathrm{P}_2\pi .
\end{equation*}
To check whether (\ref{eqn:KR}) defines a meaningful quantity, we have to check whether $\Gamma(\mu)$ is non-empty.

We shall need the following definition.

\begin{definition}\label{defin:prod}
Let $\sigma$ be an $\mathbb{R}^m$-valued Borel measure on $X$ and let $\theta$ be a Borel signed measure on $X$. A unique Borel $\mathbb{R}^m$-valued measure $\sigma\otimes \theta $ such that
\begin{equation*}
\langle \sigma\otimes\theta,v\rangle= \langle \sigma, v\rangle \otimes \theta
\end{equation*}
for all $v\in\mathbb{R}^m$ we shall call a \emph{product measure}.
Here  $\langle \sigma,v \rangle \otimes \theta$ is the usual product measure of $\mathbb{R}$-valued measures.
\end{definition}

\begin{remark}
It is clear that the product measure exists. Analogously we define the product measure $\theta\otimes \sigma$ for a Borel signed measure $\sigma$ and a Borel $\mathbb{R}^m$-valued measure $\theta$.
\end{remark}

\begin{proposition}
$\Gamma(\mu)$ is non-empty if and only if 
\begin{equation}\label{eqn:equam}
\mu(X)=0.
\end{equation}
\end{proposition}
\begin{proof}
Clearly, if there exists $\pi\in\Gamma(\mu)$, then
\begin{equation*}
\mu(X)=\mathrm{P}_1\pi(X)-\mathrm{P}_2\pi(X)=\pi(X\times X)-\pi(X\times X)=0,
\end{equation*}
so the condition (\ref{eqn:equam}) is satisfied. Conversely, assume that (\ref{eqn:equam}) holds true.
Let $\nu$ be any Borel probability measure on $X$. Set
\begin{equation*}
\pi=\mu\otimes \nu.
\end{equation*}
Here $\mu\otimes \nu$ is the product measure; see Definition~\ref{defin:prod}. 
Then for any Borel set $A\subset X$ and any vector $v\in\mathbb{R}^m$, we have 
\begin{equation*}
\langle \pi(A\times X)-\pi(X\times A),v\rangle=\langle \mu(A),v\rangle -\langle \mu(X),v\rangle \nu(A)=\langle\mu(A),v\rangle.
\end{equation*}
This is to say, $\mathrm{P}_1\pi-\mathrm{P}_2\pi=\mu$.
\qed
\end{proof}

The quantity defined by (\ref{eqn:KR}) we shall call the Kantorovich--Rubinstein norm of $\mu$.

\begin{proposition}
Assume that $\mu(X)=0$. Then $\mathcal{I}(\mu)<\infty$ provided that
\begin{equation}\label{eqn:moment}
\int_{\mathbb{R}^n}d(x,x_0)d\norm{\mu}(x)<\infty
\end{equation}
for some (equivalently: any) $x_0\in X$.
\end{proposition}
\begin{proof}
Define
\begin{equation*}
\pi=\mu\otimes \delta_{x_0}.
\end{equation*}
Here $\delta_{x_0}$ is a probability measure such that $\delta_{x_0}(\{x_0\})=1$.
Then $\pi\in\Gamma(\mu)$ and
\begin{equation}\label{eqn:finite}
\int_{X\times X}d(x,y)d\norm{\pi}(x,y)\leq \int_{X}d(x,x_0)d\norm{\mu}(x).
\end{equation}
This shows that $\mathcal{I}(\mu)<\infty$, provided that (\ref{eqn:moment}) is satisfied.
The equivalence of finiteness of
\begin{equation*}
\int_{\mathbb{R}^n}d(x,y)d\norm{\mu}(x)<\infty
\end{equation*}
for any $y\in X$ follows by the triangle inequality. \qed
\end{proof}

\begin{definition}
We define the \emph{Wasserstein space} $\mathcal{W}(X,\mathbb{R}^m)$ of all Borel measures $\mu$ on $X$ with values in $\mathbb{R}^m$ such that 
\begin{equation*}
\mu(X)=0 \text{ and }
\int_{X}d(x,x_0)d\norm{\mu}(x)<\infty
\end{equation*}
for some $x_0\in X$.
We endow it with a norm $\norm{\mu}_{\mathcal{W}(X,\mathbb{R}^m)}=\mathcal{I}(\mu)$.
\end{definition}

Before we proceed let us recall the following definition.

We say that a non-negative Borel measure $\mu$ on $X$ is \emph{inner regular} if for any Borel set $B\subset X$ we have
\begin{equation*}
\mu(B)=\sup\{\mu(K)\mid K\subset B, K \text{ is a compact set}\}.
\end{equation*}

Let us note that Ulam's lemma tells that any finite Borel measure on a Polish space is inner regular.

\begin{lemma}\label{lem:lips}
Suppose that $X$ is a Polish space. Let $\mu$ be a $\mathbb{R}^m$-valued Borel measure in $\mathcal{W}(X,\mathbb{R}^m)$. Suppose that for any Lipschitz function $u\colon X\to\mathbb{R}^m$
\begin{equation*}
\int_{X}\langle u,  d\mu\rangle =0.
\end{equation*}
Then $\mu=0$.
\end{lemma}
\begin{proof}
We may assume that $m=1$. Let $\mu=\mu_+-\mu_-$ be the Hahn--Jordan decomposition of $\mu$. There exists two disjoint, Borel sets $A,B\subset X$ with $\mu_+(A^c)=0$ and $\mu_-(B^c)=0$. Choose any Borel set $E\subset A$. As any finite measure on $X$ is inner regular, for any $\epsilon>0$, there exists a compact set $K\subset E$ such that
\begin{equation*}
\mu_+(E)\leq \mu_+(K)+\epsilon.
\end{equation*}
Define a function $u_{\epsilon}$ by the formula
\begin{equation*}
u_{\epsilon}(x)=(1-\frac{1}{\epsilon}\mathrm{dist}(x,K))\vee 0.
\end{equation*}
Then $u_{\epsilon}$ is Lipschitz, equal to one on $K$ and equal to zero on the complement of \begin{equation*}
K_{\epsilon}=\{x\in X\mid\mathrm{dist}(x,K)\leq \epsilon\}.
\end{equation*}
Thus
\begin{equation*}
0=\int_{X}u_{\epsilon}d\mu=\mu_+(K)+\int_{K_{\epsilon}\setminus K}u_{\epsilon}d\mu,
\end{equation*}
Therefore, by the above,
\begin{equation*}
\mu_+(E)\leq \epsilon+\mu_+(K)\leq \epsilon+\mu_-(K_{\epsilon}\setminus K).
\end{equation*}
Letting $\epsilon$ tend to zero, we get $\mu_+(E)=0$. It follows that $\mu_+=0$. Analogously, $\mu_-=0$. This is to say, $\mu=0$.\qed
\end{proof}

\begin{remark}
In what follows, we shall always assume that underlying space $X$ is a Polish space.
\end{remark}

\begin{proposition}
The function $\mathcal{W}(X,\mathbb{R}^m)\ni\mu\mapsto\norm{\mu}_{\mathcal{W}(X,\mathbb{R}^m)}\in\mathbb{R}$ is a norm.
\end{proposition}
\begin{proof}
Let us first check that 
\begin{equation}\label{eqn:nondeg}
\norm{\mu}_{\mathcal{W}(X,\mathbb{R}^m)}=0 \text{ if and only if } \mu=0.
\end{equation}
If $\mu=0$, then $\pi=0$ belongs to $\Gamma(\mu)$, so $\norm{\mu}_{\mathcal{W}(X,\mathbb{R}^m)}=0$. Conversely, assume that $\norm{\mu}_{\mathcal{W}(X,\mathbb{R}^m)}=0$.
Choose any $L$-Lipschitz function 
\begin{equation*}
u\colon X\to\mathbb{R}^m.
\end{equation*}
Then for any $\pi\in\Gamma(\mu)$ we have
\begin{equation*}
\Big\lvert\int_{X}\langle u, d\mu\rangle\Big\rvert= \Big\lvert\int_{X\times X}\langle u(x)-u(y), d\pi(x,y)\rangle\Big\rvert\leq L \int_{X\times X}d(x,y)d\norm{\pi}(x,y).
\end{equation*}
Therefore if $\norm{\mu}_{\mathcal{W}(X,\mathbb{R}^m)}=0$, then
\begin{equation*}
\int_{X}\langle u, d\mu\rangle=0.
\end{equation*}
It follows by Lemma \ref{lem:lips}, that $\mu=0$.
Homogeneity of  $\norm{\cdot}_{\mathcal{W}(X,\mathbb{R}^m)}$ is clear. Let us show that the triangle inequality holds. 
For this choose measures $\mu,\nu\in \mathcal{W}(X,\mathbb{R}^m)$ and any measures $\pi\in \Gamma(\mu)$ and $\rho\in\Gamma(\nu)$. Then
\begin{equation*}
\mu+\nu=\mathrm{P}_1(\pi+\rho)-\mathrm{P}_2(\pi+\rho),
\end{equation*}
so that $\pi+\rho\in\Gamma(\mu+\nu)$.
It follows that 
\begin{equation*}
\begin{aligned}
\norm{\mu+\nu}_{\mathcal{W}(X,\mathbb{R}^m)}&\leq \int_{X\times X}d(x,y)d\norm{\pi+\rho}(x,y)\leq\\
&\leq \int_{X\times X}d(x,y)d\norm{\pi}(x,y)+ \int_{X\times X}d(x,y)d\norm{\rho}(x,y).
\end{aligned}
\end{equation*}
Taking infimum over all $\pi,\rho$ we see that the triangle inequality holds true.\qed
\end{proof}

\begin{proposition}\label{pro:density}
The linear space $\mathcal{U}$ of measures of the form
\begin{equation*}
\sum_{i=1}^n \delta_{x_i} v_i
\end{equation*}
for $x_i\in X$ and $v_i\in\mathbb{R}^m$, $i=1,\dotsc,n$, such that $\sum_{i=1}^n v_i=0$, is dense in $\mathcal{W}(X,\mathbb{R}^m)$.
\end{proposition}
\begin{proof}
Choose any measure $\mu\in\mathcal{W}(X,\mathbb{R}^m)$. Choose any $\epsilon>0$. Choose any point $x_0\in X$ and a compact set $K$ such that
\begin{equation*}
\int_{K^c}d(x,x_0)d\norm{\mu}(x)\leq \epsilon.
\end{equation*}
Choose pairwise disjoint Borel sets $A_1,A_2,\dotsc,A_k \subset K$ such that the diameter of each is at most $\epsilon$ and
\begin{equation*}
K=\bigcup_{i=1}^kA_i.
\end{equation*}
Consider the restrictions $\mu_i=\mu|_{A_i}$ of the measure $\mu$ to the sets $A_i$, $i=1,2,\dotsc,k$. Choose any points $x_i\in A_i$. Then, as 
\begin{equation*}
\pi_i=\mu_i\otimes\delta_{x_i}\in\Gamma(\mu_i-\mu_i(X)\delta_{x_i}),
\end{equation*}
we have
\begin{equation*}
\norm{\mu_i- \mu_i(X) \delta_{x_i}}_{\mathcal{W}(X,\mathbb{R}^m)}\leq \int_{X}d(y,x_i)d\norm{\mu_i}(y)\leq \epsilon \norm{\mu}(A_i).
\end{equation*}
Let $A_0=K^c$ and let $\mu_0=\mu|_{A_0}$. Then 
\begin{equation*}
\pi_0=\mu_0\otimes \delta_{x_0}\in\Gamma(\mu_0-\mu_0(X)\delta_{x_0}),
\end{equation*}
so
\begin{equation*}
\norm{\mu_0-\mu_0(X)\delta_{x_0}}_{\mathcal{W}(X,\mathbb{R}^m)}\leq \int_{X}d(x,x_0)d\norm{\mu_0}(x)\leq \epsilon.
\end{equation*}
Set 
\begin{equation*}
\nu=\sum_{i=0}^k\mu(A_i)\delta_{x_i}.
\end{equation*}
Then $\nu\in \mathcal{U}$. By the triangle inequality
\begin{equation*}
\begin{aligned}
&\norm{\mu-\nu}_{\mathcal{W}(X,\mathbb{R}^m)}\leq \sum_{i=0}^k \norm{\mu_i-\mu_i(X)\delta_{x_i}}_{\mathcal{W}(X,\mathbb{R}^m)}\leq\\
&\leq \epsilon \sum_{i=1}^k\norm{\mu(A_i)}+\epsilon \leq \epsilon ( \norm{\mu}(X)+1).
\end{aligned}
\end{equation*}
This concludes the proof.
\qed
\end{proof}

\begin{corollary}
If $X$ is separable, then so is the Wasserstein space $\mathcal{W}(X,\mathbb{R}^m)$.
\end{corollary}
\begin{proof}
Fix $n\in\mathbb{N}$. Choose a countable dense subset $A\subset X$ and a set
\begin{equation}\label{eqn:dense}
B\subset \Big\{(w_1,\dotsc,w_n)\in\mathbb{R}^m\times\dotsc\mathbb{R}^m\mid \sum_{i=1}^nw_i=0\Big\}
\end{equation}
which is countable and dense in the set on the right-hand side of (\ref{eqn:dense}).
Consider a measure $\mu$ given by
\begin{equation*}
\mu=\sum_{i=1}^n \delta_{x_i} v_i
\end{equation*}
for $x_i\in X$ and $v_i\in\mathbb{R}^m$, $i=1,\dotsc,n$, such that $\sum_{i=1}^n v_i=0$. Choose $\epsilon>0$ and $\tilde{x}_i\in A$, $i=1,\dotsc,n$, and $(\tilde{v}_i)_{i=1}^n\in B$, such that for $i=1,\dotsc,n$
\begin{equation*}
d(x_i,\tilde{x}_i)<\epsilon \text{ and }\norm{v_i-\tilde{v}_i}<\epsilon \text{ and } \sum_{i=1}^n\tilde{v}_i=0.
\end{equation*} 
Set 
\begin{equation*}
\tilde{\mu}=\sum_{i=1}^n \delta_{\tilde{x}_i} \tilde{v}_i.
\end{equation*}
Then
\begin{equation*}
\norm{\mu-\tilde{\mu}}_{\mathcal{W}(X,\mathbb{R}^m)}\leq \Big\lVert \sum_{i=1}^n \delta_{x_i}(v_i-\tilde{v}_i)\Big\rVert_{\mathcal{W}(X,\mathbb{R}^m)}+\Big\lVert \sum_{i=1}^n (\delta_{x_i}-\delta_{\tilde{x}_i})v_i\Big\rVert_{\mathcal{W}(X,\mathbb{R}^m)}  
\end{equation*}
Choose any $x_0\in X$. Taking 
\begin{equation*}
\pi=\sum_{i=1}^n \delta_{x_i} \otimes \delta_{x_0}(v_i-\tilde{v}_i)\text{ and }
\rho=\sum_{i=1}^n (\delta_{x_i}\otimes \delta_{\tilde{x}_i})v_i
\end{equation*}
we see that
\begin{equation*}
\Big\lVert \sum_{i=1}^n \delta_{x_i}(v_i-\tilde{v}_i)\Big\rVert_{\mathcal{W}(X,\mathbb{R}^m)}\leq \epsilon \sum_{i=1}^nd(x_i,x_0)
\end{equation*}
and
\begin{equation*}
\Big\lVert \sum_{i=1}^n (\delta_{x_i}-\delta_{\tilde{x}_i})v_i\Big\rVert_{\mathcal{W}(X,\mathbb{R}^m)}\leq \epsilon \sum_{i=1}^n\norm{v_i}.
\end{equation*}
The conclusion follows now from Proposition \ref{pro:density}.\qed
\end{proof}

\begin{definition}
Choose any $x_0\in X$. Define
\begin{equation*}
\mathcal{L}(X,\mathbb{R}^m)=\{u\colon X\to\mathbb{R}^m\mid u \text{ is Lipschitz and } u(x_0)=0\},
\end{equation*}
i.e. the Banach space of $\mathbb{R}^m$-valued Lipschitz functions on $X$ taking value zero at $x_0$,
with norm 
\begin{equation*}
\norm{u}_{\mathcal{L}(X,\mathbb{R}^m)}=\sup\bigg\{\frac{\norm{u(x)-u(y)}}{d(x,y)}\mid x,y\in X, x\neq y\bigg\}.
\end{equation*}
\end{definition}

\begin{theorem}\label{thm:reciprocals}
Define 
\begin{equation*}
T\colon\mathcal{L}(X,\mathbb{R}^m)\to\mathcal{W}(X,\mathbb{R}^m)^*
\end{equation*}
and 
\begin{equation*}
S\colon\mathcal{W}(X,\mathbb{R}^m)^*\to \mathcal{L}(X,\mathbb{R}^m)
\end{equation*}
by 
\begin{equation}\label{eqn:isodual}
T(u)(\mu)=\int_{X}\langle u,d\mu\rangle
\end{equation}
and
\begin{equation}\label{eqn:dualiso}
\langle S(\lambda)(x),w\rangle=\lambda((\delta_x-\delta_{x_0})w),
\end{equation}
for any $w\in\mathbb{R}^m$. Then $S, T$ are mutual reciprocals and establish an isometric isomorphism of $\mathcal{L}(X,\mathbb{R}^m)$ and $\mathcal{W}(X,\mathbb{R}^m)^*$. 
\end{theorem}
\begin{proof}
Choose any $\pi\in\Gamma(\mu)$. Then $\mathrm{P}_1\pi-\mathrm{P}_2\pi=\mu$. Thus, if $u$ is a Lipschitz map, then
\begin{equation*}
\bigg\lvert \int_X\langle u, d\mu \rangle \bigg\rvert=\bigg\lvert \int_X\langle u(x)-u(y), d\pi(x,y) \rangle \bigg\rvert\leq \norm{u}_{\mathcal{L}(X,\mathbb{R}^m)}\int_X d(x,y) d\norm{\pi}(x,y).
\end{equation*}
Taking infimum over all $\pi\in\Gamma(\mu)$, we see that 
\begin{equation*}
\bigg\lvert \int_X\langle u, d\mu \rangle \bigg\rvert\leq \norm{u}_{\mathcal{L}(X,\mathbb{R}^m)}\norm{\mu}_{\mathcal{W}(X,\mathbb{R}^m)}.
\end{equation*}
The above calculation shows that the formula (\ref{eqn:isodual}) defines a continuous functional of norm at most $\norm{u}_{\mathcal{L}(X,\mathbb{R}^m)}$. If $w\in\mathbb{R}^m$ if of norm one and $x,y\in X$, $x\neq y$, then for 
\begin{equation}\label{eqn:mu}
\mu_{x,y,w}=\frac{\delta_x-\delta_y}{d(x,y)}w
\end{equation}
we have $\norm{\mu_{x,y,w}}_{\mathcal{W}(X,\mathbb{R}^m)}\leq 1$ and for any $u\in\mathcal{L}(X,\mathbb{R}^m)$
\begin{equation*}
\int_{\mathbb{R}^n}\langle u,d\mu_{x,y,w}\rangle=\frac{\langle w, u(x)-u(y)\rangle}{d(x,y)}.
\end{equation*}
Thus
\begin{equation*}
\norm{u}_{\mathcal{L}(X,\mathbb{R}^m)}=\norm{T(u)}.
\end{equation*}
We shall now show that $T\circ S=\mathrm{Id}$. Take any functional $\lambda\in\mathcal{W}(X,\mathbb{R}^m)^*$. Set 
\begin{equation*}
\sigma_{x,w}=(\delta_x-\delta_{x_0}) w.
\end{equation*}
Then  $S(\lambda)\colon X\to\mathbb{R}^m$ is defined by the formula
\begin{equation*}
\langle S(\lambda)(x),w\rangle=\lambda(\sigma_{x,w}).
\end{equation*}
It is clear that the above formula defines $S(\lambda)$ uniquely. Then we claim that map $v=S(\lambda)$ is $\norm{\lambda}$-Lipschitz. Indeed
\begin{equation*}
\norm{v(x)-v(y)}=\sup\{\langle v(x)-v(y),w\rangle \mid w\in\mathbb{R}^m, \norm{w}=1\},
\end{equation*}
and as
\begin{equation*}
\langle v(x)-v(y),w\rangle = \lambda (\sigma_{x,w}-\sigma_{y,w})\leq \norm{\lambda}\norm{\sigma_{x,w}-\sigma_{y,w}}_{\mathcal{W}(X,\mathbb{R}^m)}
\end{equation*}
we see that 
\begin{equation*}
\norm{v(x)-v(y)}\leq\norm{\lambda}d(x,y), \text{ since } \norm{\sigma_{x,w}-\sigma_{y,w}}_{\mathcal{W}(X,\mathbb{R}^m)}\leq d(x,y).
\end{equation*}
Suppose that $\nu=(\delta_x-\delta_y)z$. We compute
\begin{equation*}
T(v)(\nu)=\int_X\langle v,d\nu\rangle=\int_X\langle v,z\rangle  d(\delta_x-\delta_y)= \lambda (\sigma_{x,z}-\sigma_{y,z})=\lambda(\nu).
\end{equation*}
We see that $T(S(\lambda))$ and $\lambda$ are equal on the set spanned by $(\delta_x-\delta_y)z$, where $x,y\in X$, $z\in\mathbb{R}^m$. By Proposition \ref{pro:density}, we see that $T(S(\lambda))$ and $\lambda$ are equal on $\mathcal{W}(X,\mathbb{R}^m)$. 

Let us show also that $S\circ T=\mathrm{Id}$. Choose any $w\in\mathbb{R}^m$ and any map $u\in\mathcal{L}(X,\mathbb{R}^m)$. Then
\begin{equation*}
\langle S(T(u))(x),w\rangle = T(u)((\delta_x-\delta_{x_0})w)=\int_X\langle u, d(\delta_x-\delta_{x_0})w\rangle=\langle u(x),w\rangle,
\end{equation*}
as $u(x_0)=0$. Therefore $S(T(u))=u$.\qed
\end{proof}

\begin{theorem}\label{thm:lip}
For any $\mu\in \mathcal{W}(X,\mathbb{R}^m)$
\begin{equation}\label{eqn:Kantorovich-Rubinsetin}
\sup\bigg\{\int_X\langle u, d\mu \rangle \mid u\colon X\to\mathbb{R}^m \text{ is } 1\text{-Lipschitz}\bigg \}=\norm{\mu}_{\mathcal{W}(X,\mathbb{R}^m)}.
\end{equation}
Moreover, there exists $1$-Lipschitz function $u_0$ such that
\begin{equation}\label{eqn:maxip}
\sup\bigg\{\int_X\langle u, d\mu \rangle \mid u\colon X\to\mathbb{R}^m \text{ is } 1\text{-Lipschitz}\bigg \}=\int_X\langle u_0, d\mu\rangle.
\end{equation}
\end{theorem}
\begin{proof}
Notice first that the left-hand side of (\ref{eqn:Kantorovich-Rubinsetin}) is clearly at most the right-hand side of (\ref{eqn:Kantorovich-Rubinsetin}). Take any $\mu\in\mathcal{W}(X,\mathbb{R}^m)$. Then by the Hahn--Banach theorem there exists a continuous linear functional $\lambda$ of norm one such that 
\begin{equation*}
\lambda(\mu)=\norm{\mu}_{\mathcal{W}(X,\mathbb{R}^m)}.
\end{equation*}
By Theorem \ref{thm:reciprocals}, we know that $\lambda$ is of the form
\begin{equation*}
\lambda(\mu)=\int_X\langle u_0,d\mu\rangle
\end{equation*}
for some Lipschitz map $u_0$. The Lipschitz constant of $u_0$ is equal to one, as
\begin{equation*}\norm{u_0}_{\mathcal{L}(X,\mathbb{R}^m)}=\norm{\lambda}=1.
\end{equation*}
This completes the proof.\qed
\end{proof}

\begin{definition}
Any $1$-Lipschitz function $u\colon X\to\mathbb{R}^m$ such that (\ref{eqn:maxip}) holds we shall call an \emph{optimal potential} of measure $\mu$.
\end{definition}

\begin{definition}
A measure $\pi\in \Gamma(\mu)$ such that
\begin{equation*}
\norm{\mu}_{\mathcal{W}(X,\mathbb{R}^m)}=\int_{X\times X}d(x,y)d\norm{\pi}(x,y)
\end{equation*}
we shall call an \emph{optimal transport} for $\mu$.
\end{definition}

\begin{theorem}\label{thm:ae}
Let $\mu\in\mathcal{W}(X,\mathbb{R}^m)$. Let $u\in\mathcal{L}(X,\mathbb{R}^m) $ be a $1$-Lipschitz map. Let $\pi\in\Gamma(\mu)$. The following conditions are equivalent:
\begin{enumerate}[i)]
\item\label{i:opti}
\begin{equation*}
\int_X \langle u,d\mu\rangle=\int_{X\times X}d(x,y)d\norm{\pi}(x,y)=\norm{\mu}_{\mathcal{W}(X,\mathbb{R}^m)} ,
\end{equation*}
\item\label{i:borel}
\begin{equation*}
\int_A \langle u(x)-u(y),d\pi(x,y) \rangle=\int_A d(x,y)d\norm{\pi}(x,y)
\end{equation*}
for any Borel set $A\subset X\times X$,
\item \label{i:eren}
\begin{equation*}
\int_X \langle u,d\mu\rangle=\int_{X\times X} d(x,y)d\norm{\pi}(x,y),
\end{equation*}
\item\label{i:optopt}
$u$ is an optimal potential for $\mu$ and $\pi$ is an optimal transport for $\mu$.
\end{enumerate}
Moreover, if the above conditions hold, then
\begin{equation*}
\norm{u(x)-u(y)}=d(x,y)
\end{equation*}
$\norm{\pi}$-almost everywhere.
\end{theorem}
\begin{proof}
Assume that \ref{i:eren}) holds. Observe that
\begin{equation*}
\int_X\langle u, d\mu\rangle=\int_{X\times X} \langle u(x)-u(y),d\pi(x,y) \rangle.
\end{equation*}
As
\begin{equation*}
\int_X\langle u, d\mu\rangle\leq \norm{\mu}_{\mathcal{W}(X,\mathbb{R}^m)}\leq \int_{X\times X}d(x,y)d\norm{\pi}(x,y),
\end{equation*}
then by \ref{i:eren}) we see that in the above inequalities we have equalities. This is to say, \ref{i:opti}) holds true. 

Suppose now that \ref{i:opti}) holds. Clearly
\begin{equation*}
\int_A \langle u(x)-u(y),d\pi(x,y) \rangle\leq\int_A d(x,y)d\norm{\pi}(x,y).
\end{equation*}
If we had strict inequality in \ref{i:borel}) for some Borel set $A\subset X\times X$, then the above computations show that we would get strict inequality in \ref{i:opti}).
Condition \ref{i:optopt}) is reformulation of \ref{i:opti}).
The last part of the theorem follows readily from \ref{i:borel}).\qed
\end{proof}

We say that a measure $\mu\in\mathcal{M}(Z,\mathbb{R}^m)$ is \emph{concentrated} on a subset $X\subset Z$ if there is $\norm{\mu}(Z\setminus X)=0$.

\begin{proposition}\label{pro:concent}
Assume that $\mathbb{R}^n,\mathbb{R}^m$ are equipped with Euclidean norms. Let $\mu\in \mathcal{W}(\mathbb{R}^n,\mathbb{R}^m)$ be concentrated on a set $X\subset\mathbb{R}^n$. Then
\begin{equation*}
\norm{\mu}_{\mathcal{W}(\mathbb{R}^n,\mathbb{R}^m)}=\norm{\mu}_{\mathcal{W}(X,\mathbb{R}^m)}.
\end{equation*}
\end{proposition}
\begin{proof}
The assertion is that
\begin{equation*}
\sup\Big\{\int_{\mathbb{R}^n}\langle u,d\mu\rangle\mid u\colon\mathbb{R}^n\to\mathbb{R}^m\text{ is }1\text{-Lipschitz}\Big\}
\end{equation*}
is equal to
\begin{equation*}
\sup\Big\{\int_X\langle u,d\mu\rangle\mid u\colon X\to\mathbb{R}^m\text{ is }1\text{-Lipschitz}\Big\}.
\end{equation*}
By the Kirszbraun theorem (see e.g. \cite{Kirszbraun}) any $1$-Lipschitz function $u\colon X\to\mathbb{R}^m$ extends to a $1$-Lipschitz function $\tilde{u}\colon\mathbb{R}^n\to\mathbb{R}^m$. Clearly, for any such extension
\begin{equation*}
\int_{\mathbb{R}^n}\langle \tilde{u},d\mu\rangle=\int_X\langle u,d\mu\rangle.
\end{equation*}
The assertion follows.\qed
\end{proof}

\section{Mass balance condition}\label{sec:counter}

Let us first provide an affirmative answer to the conjecture of Klartag, under the provision of the existence of optimal transport with absolutely continuous marginals of its total variation.

\begin{definition}
A leaf $\mathcal{S}$ of a $1$-Lipschitz map $u\colon\mathbb{R}^n\to\mathbb{R}^m$ is a maximal set, with respect to the order induced by inclusion, such that the restriction $u|_{\mathcal{S}}$ is an isometry. This is to say, $\mathcal{S}$ is a leaf, whenever for any $x,y\in\mathcal{S}$ there is 
\begin{equation*}
\norm{u(x)-u(y)}=\norm{x-y}
\end{equation*}
and for any $z\notin\mathcal{S}$ there exists $x\in\mathcal{S}$ such that
\begin{equation*}
\norm{u(x)-u(z)}<\norm{x-z}.
\end{equation*}
\end{definition}

It is proven in \cite{Ciosmak2} that leaves of a map $u$ that is $1$-Lipschitz with respect to Euclidean norms are closed and convex sets. Two distinct leaves may intersect at most by their relative boundaries.

\begin{definition}
Let $u\colon\mathbb{R}^n\to\mathbb{R}^m$ be a $1$-Lipschitz map of Euclidean spaces. We say that a Borel set $A\subset\mathbb{R}^n$ is a \emph{transport set} associated with $u$ if it enjoys the following property: if $x\in A$ is contained in a unique leaf of $u$ and $y\in\mathbb{R}^n$ is such that 
\begin{equation*}
\norm{u(x)-u(y)}=\norm{x-y},
\end{equation*}
then $y\in A$.
\end{definition}

Let us remark that a Borel set $A\subset\mathbb{R}^n$ that is a union of leaves of $u$ is a transport set. 

We shall denote by $B(u)$ the set of all points $x\in\mathbb{R}^n$ such that there exist at least two distinct leaves $\mathcal{S}_1,\mathcal{S}_2$ of $u$ such that $x\in\mathcal{S}_1\cap\mathcal{S}_2$. In \cite[Corollary 2.15]{Ciosmak2} it is proven that $B(u)$ is of Lebesgue measure zero.

Suppose that $\mu\in\mathcal{W}(\mathbb{R}^n,\mathbb{R}^m)$.
The following theorem shows that if there exists an optimal transport for $\mu$ such that its total variation has absolutely continuous marginals, then the conjecture of Klartag holds true. Note that such existence is clear for $m=1$, whenever $\mu$ is absolutely continuous with respect to the Lebesgue measure $\lambda$.

\begin{theorem}\label{thm:balance}
Assume that $\mathbb{R}^n,\mathbb{R}^m$ are equipped with Euclidean norms. Suppose that $\mu\in\mathcal{W}(\mathbb{R}^n,\mathbb{R}^m)$. Let $u$ be an optimal potential for $\mu$. Suppose that there exists an optimal transport $\pi$ of $\mu$ such that
\begin{equation}\label{eqn:abs}
\mathrm{P}_1\norm{\pi}\ll, \mathrm{P}_2\norm{\pi} \ll\lambda.
\end{equation}
Then for any transport set $A$ associated with $u$:
\begin{enumerate}[i)]
\item\label{i:mass}  $\mu(A)=0$,
\item\label{i:transport}
$\pi|_{A\times A}\in \Gamma(\mu|_A)$ is an optimal transport of $\mu|_A$
\item\label{i:potential}
$u$ is an optimal potential of $\mu|_A$.
\end{enumerate} 
\end{theorem}
\begin{proof}
By \cite[Corollary 2.15]{Ciosmak2} it follows that
\begin{equation*}
\lambda(B(u))=0.
\end{equation*}
Suppose that (\ref{eqn:abs}) holds true. Then
\begin{equation*}
\norm{\pi}\big(B(u)\times \mathbb{R}^n \big)=0\text{ and }\norm{\pi}\big( \mathbb{R}^n\times B(u) \big)=0.
\end{equation*}
Let 
\begin{equation*}
I=\big\{(x,y)\in\mathbb{R}^n\times\mathbb{R}^n\mid \norm{u(x)-u(y)}=\norm{x-y}\big\}.
\end{equation*}
By Theorem \ref{thm:ae}, $\norm{\pi}(I^c)=0$. 
Thus $\pi$ is concentrated on the set 
\begin{equation*}
C=I\cap\big( B(u)^c\times  B(u)^c\big).
\end{equation*}
Suppose that $(x,y)\in C$.
Then, as $A$ is a transport set, by the definition of $B(u)$,
\begin{equation}\label{eqn:equiv}
x\in A \text{ if and only if }y\in A.
\end{equation}
Let $\eta=\pi|_{A\times A}$. 
To prove \ref{i:transport}), it is enough to show that $\eta$ is an optimal transport and that  
\begin{equation*}
\eta\in\Gamma(\mu|_{A}).
\end{equation*}
For this, let $D\subset \mathbb{R}^n$ be any Borel set. Using the fact that $\pi\in\Gamma(\mu)$ and the fact that $\norm{\pi}(C^c)=0$ and (\ref{eqn:equiv}),  we have
\begin{equation*}
\begin{aligned}
&\mu(A\cap D)=\int_{\mathbb{R}^n\times\mathbb{R}^n}\Big(\mathbf{1}_{A\cap D}(x)-\mathbf{1}_{A\cap D}(y)\Big)d\pi(x,y)=\\
&=\int_{\mathbb{R}^n\times\mathbb{R}^n}\mathbf{1}_{A\times A}(x,y)\Big(\mathbf{1}_D(x)-\mathbf{1}_D(y)\Big)d\pi(x,y)=\\
&=\int_{\mathbb{R}^n\times\mathbb{R}^n}\Big(\mathbf{1}_D(x)-\mathbf{1}_D(y)\Big)d\eta(x,y)=\mathrm{P}_1\eta(D)-\mathrm{P}_2\eta(D).
\end{aligned}
\end{equation*}
It follows that $\pi|_{A\times  A}\in\Gamma(\mu|_A)$.
Then
\begin{equation}\label{eqn:comput}
\int_{A}\langle u,d\mu\rangle =\int_{\mathbb{R}^n\times\mathbb{R}^n} \mathbf{1}_C(x,y)\Big\langle \mathbf{1}_A(x) u(x)-\mathbf{1}_A(y)u(y),d\pi(x,y)\Big\rangle .
\end{equation}
Therefore, by (\ref{eqn:equiv}),
\begin{equation*}
\int_{A}\langle u,d\mu\rangle
=\int_{\mathbb{R}^n\times\mathbb{R}^n}\mathbf{1}_{A\times A}(x,y)\Big\langle u(x)-u(y),d\pi(x,y)\Big\rangle .
\end{equation*}
By condition \ref{i:borel}) of Theorem \ref{thm:ae} we see that
\begin{equation*}
\int_{A}\langle u,d\mu\rangle=\int_{A\times A}\norm{x-y}d\norm{\pi}(x,y).
\end{equation*}
Theorem \ref{thm:ae}, condition \ref{i:eren}), tells us that $\pi|_{A\times A}$ is an optimal transport and $u$ is an optimal potential. Also $\mu(A)=0$, as $\pi|_{A\times A}\in\Gamma(\mu|_A)$. This completes the proof.
\qed
\end{proof}

We shall now provide necessary tools for the aforementioned counterexample to the conjecture of Klartag.

In fact we shall provide a more general theorem for which we shall consider locally uniformly closed subsets subsets $\mathcal{F}$ of $1$-Lipschitz maps of $\mathbb{R}^n$ to $\mathbb{R}^m$ endowed with norms which are not necessarily Euclidean. Suppose that a measure $\mu$ belongs to $\mathcal{W}(\mathbb{R}^n,\mathbb{R}^m)$. We consider supremum of integrals
\begin{equation}\label{eqn:maxil}
\int_{\mathbb{R}^n}\langle u,d\mu\rangle
\end{equation}
taken over all $u\in\mathcal{F}$.
An optimal $u_0\in\mathcal{F}$, i.e. the map that satisfies 
\begin{equation*}
\int_{\mathbb{R}^n}\langle u_0,d\mu\rangle=\sup\Big\{\int_{\mathbb{R}^n}\langle u,d\mu\rangle\mid u\in\mathcal{F}\Big\},
\end{equation*}
we shall call an $\mathcal{F}$-optimal potential of $\mu$.

\begin{lemma}\label{lem:attain}
Let $X\subset\mathbb{R}^n$ be a compact set. Suppose that $(\mu_k)_{k=1}^{\infty}\subset\mathcal{W}(\mathbb{R}^n,\mathbb{R}^m)$ are all supported on $X$ and converge weakly* to $\mu_0\in\mathcal{W}(\mathbb{R}^n,\mathbb{R}^m)$, i.e. for any 
continuous and bounded function $g\colon \mathbb{R}^n\to\mathbb{R}^m$ we have
\begin{equation*}
\lim_{k\to\infty}\int_{\mathbb{R}^n}\langle g,d\mu_k\rangle=\int_{\mathbb{R}^n}\langle g,d\mu_0\rangle.
\end{equation*}
Suppose that for $k=1,2,\dotsc,$ $u_k\in\mathcal{F}$ is an $\mathcal{F}$-optimal potential of $\mu_k$
and that $u_k$ converge locally uniformly to $u_0\colon \mathbb{R}^n\to\mathbb{R}^m$. Then $u_0$ is an $\mathcal{F}$-optimal potential of $\mu_0$.
\end{lemma}
\begin{proof}
By the assumption, for any continuous and bounded map $g\colon \mathbb{R}^n\to\mathbb{R}^m$, we have
\begin{equation*}
\lim_{k\to\infty}\int_{\mathbb{R}^n}\langle g,d(\mu_k-\mu_0)\rangle=0.
\end{equation*}
In particular, as $\mu_k$ are all supported on $X$, we have
\begin{equation*}
\lim_{k\to\infty}\int_{\mathbb{R}^n}\langle u_0,d(\mu_k-\mu_0)\rangle=0.
\end{equation*}
By the Banach--Steinhaus theorem, the sequence $(\mu_k)_{k=1}^{\infty}$ is bounded in the total variation norm.
Hence, by uniform convergence on $X$,
\begin{equation*}
\lim_{k\to\infty}\int_{\mathbb{R}^n}\langle u_k-u_0,d\mu_k\rangle=0.
\end{equation*}
It follows that 
\begin{equation*}
\int_{\mathbb{R}^n}\langle u_k,d\mu_k\rangle=\int_{\mathbb{R}^n}\langle u_0,d\mu_k\rangle+\int_{\mathbb{R}^n}\langle u_k-u_0,d\mu_k\rangle
\end{equation*}
converges to $\int_{\mathbb{R}^n}\langle u_0,d\mu_0\rangle$.
As for any $1$-Lipschitz map $h\in\mathcal{F}$ we have
\begin{equation*}
\int_{\mathbb{R}^n}\langle h,d\mu_k\rangle\leq \int_{\mathbb{R}^n}\langle u_k,d\mu_k\rangle.
\end{equation*}
we also have
\begin{equation*}
\int_{\mathbb{R}^n}\langle h,d\mu_0\rangle\leq \int_{\mathbb{R}^n}\langle u_0,d\mu_0\rangle.
\end{equation*}
The proof is complete.\qed
\end{proof}

Below we shall denote by $B(x,\epsilon)$ an open ball of radius $\epsilon>0$ centred at $x\in\mathbb{R}^n$.

\begin{lemma}\label{lem:non-trivial}
Let $m\leq n$. Let $\mu\in\mathcal{W}(\mathbb{R}^n,\mathbb{R}^m)$ and let $u$ be an optimal potential of $\mu$. Let $A$ be the union of all leaves of dimension at least one. Then $A$ is Borel measurable. 
Suppose that there exists an optimal transport $\pi$ for $\mu$ or that any transport set of $u$ is of $\mu$-measure zero.  Then
\begin{equation*}
\norm{\mu}(A^c)=0.
\end{equation*}
\end{lemma}
\begin{proof}
Observe that
\begin{equation*}
A=\bigcup_{n=1}^{\infty}\Big\{x\in\mathbb{R}^n\mid \sup\Big\{\frac{\norm{u(x)-u(y)}}{\norm{x-y}}\mid y\in \mathrm{cl}B(x,n)\setminus B(x,1/n)\Big\}=1\Big\}.
\end{equation*}
The function
\begin{equation*}
\mathbb{R}^n\ni x\mapsto \sup\Big\{\frac{\norm{u(x)-u(y)}}{\norm{x-y}}\mid y\in \mathrm{cl}B(x,n)\setminus B(x,1/n)\Big\}\in \mathbb{R}
\end{equation*}
is lower semi-continuous, hence Borel measurable. Thus, $A$ is Borel measurable.
Suppose that there exists an optimal transport $\pi$ for $\mu$. By Theorem \ref{thm:ae}, $\pi$ is supported on the set
\begin{equation*}
I=\big\{(x,y)\in\mathbb{R}^n\times\mathbb{R}^n \mid \norm{u(x)-u(y)}=\norm{x-y}\big\}.
\end{equation*}
As $\mu=\mathrm{P}_1\pi-\mathrm{P}_2\pi$, for any Borel set $B\subset A^c$, we have
\begin{equation*}
\mu(B)=\pi(B\times\mathbb{R}^n)-\pi(\mathbb{R}^n\times B)=0,
\end{equation*}
for if $B\subset A^c$, then 
\begin{equation*}
\big(B\times\mathbb{R}^n\big)\cap I\subset \{(x,x)\mid x\in\mathbb{R}^n\}\text{ and }\big(\mathbb{R}^n\times B\big) \cap I\subset \{(x,x)\mid x\in\mathbb{R}^n\}.
\end{equation*}
Suppose now that any transport set for $u$ is of $\mu$ measure zero. Observe that any Borel set $B\subset A^c$ is a transport set. The conclusion follows.\qed
\end{proof}

In the theorem below we shall provide a counterexample to the conjecture of Klartag.

\begin{theorem}\label{thm:nonzero}
Assume that $m>1$. There exists an absolutely continuous measure $\mu\in\mathcal{W}(\mathbb{R}^n,\mathbb{R}^m)$ for which there exists a transport set associated with an optimal potential of $\mu$ with non-zero measure $\mu$.  

In particular, there is no optimal transport $\pi$ for $\mu$ such that
\begin{equation*}
\mathrm{P}_1\norm{\pi}\ll\lambda \text{ and }\mathrm{P}_2\norm{\pi}\ll\lambda.
\end{equation*}
\end{theorem}
\begin{proof}

Choose any $v_1,\dotsc,v_{m+1}\in\mathbb{R}^m$ such that
\begin{equation*}
\sum_{i=1}^{m+1}v_i=0
\end{equation*} 
and that are affinely independent.
For $\epsilon>0$ set
\begin{equation*}
\mu_{\epsilon}=\frac1{\lambda(B(0,\epsilon))}\sum_{i=1}^{m+1} \lambda|_{B(x_i,\epsilon)}v_i,
\end{equation*}
where $x_1,\dotsc,x_{m+1}\in\mathbb{R}^n$ are pairwise distinct points to be specified later. Here $\lambda$ denotes the Lebesgue measure on $\mathbb{R}^n$. 
Then $\mu_{\epsilon}\in\mathcal{W}(\mathbb{R}^n,\mathbb{R}^m)$. Suppose that for some sequence $(\epsilon_k)_{k=1}^{\infty}$ converging to zero there is
\begin{equation*}
\mu_{\epsilon_k}(C_k)=0
\end{equation*}
for any transport set $C_k$ of $u_k$, where $u_k\colon\mathbb{R}^n\to\mathbb{R}^m$ is an optimal potential of $\mu_{\epsilon_k}$.
For $k\in\mathbb{N}$ and $i=1,\dotsc,m+1$ consider the union $N_{ik}$ of all non-trivial leaves of $u_k$ -- i.e. of dimension at least one -- that intersect $\mathrm{cl}B(x_i,\epsilon_k)$. Then $N_{ik}$ is a transport set.
Its Borel measurability follows by Lemma \ref{lem:non-trivial}. Indeed, denote $B=\mathrm{cl}B(x_i,\epsilon_k)$; then the function
\begin{equation*}
\mathbb{R}^n\setminus B\ni x\mapsto \sup\Big\{\frac{\norm{u_k(x)-u_k(y)}}{\norm{x-y}}\mid y\in B\Big\}\in\mathbb{R}
\end{equation*}
is lower semi-continuous and therefore
\begin{equation*}
N_{ik}=\Big\{x\in \mathbb{R}^n\setminus B\mid\sup\Big\{\frac{\norm{u_k(x)-u_k(y)}}{\norm{x-y}}\mid y\in B\Big\}=1\Big\}\cup \big(B\cap A_k\big)
\end{equation*}
is a Borel set. Here $A_k$ is a set of all leaves of dimension at least one corresponding to $u_k$, c.f. Lemma \ref{lem:non-trivial}.
Thus $\mu_{\epsilon_k}(N_{ik})=0$. Hence, 
\begin{equation}\label{eqn:sum}
\sum_{j=1}^{m+1}v_j\lambda(B(x_j,\epsilon_k)\cap N_{ik})=0.
\end{equation}
As $\mu_{\epsilon_k}$, by Lemma \ref{lem:non-trivial}, is concentrated on non-trivial leaves of $u_k$, we have for 
\begin{equation*}
\frac{\lambda(B(x_i,\epsilon_k)\cap N_{ik})}{\lambda(B(0,\epsilon_k))}v_i=\mu_{\epsilon_k}(B(x_i,\epsilon_k)\cap N_{ik})=\mu_{\epsilon_k}(B(x_i,\epsilon_k))=v_i.
\end{equation*}
By (\ref{eqn:sum}) and assumption on the vectors $v_1,\dotsc,v_{m+1}$
\begin{equation*}
\lambda(B(x_j,\epsilon_k)\cap N_{ik})=\lambda(B(0,\epsilon_k)) \text{ for all }j=1,\dotsc,m+1.
\end{equation*}
Thus we infer that for any $k\in\mathbb{N}$ and for all $r,s=1,\dotsc,m+1$, $r\neq s$, there exist points
\begin{equation*}
(x_{rs}^k,x_{sr}^k)\in B(x_r,\epsilon_k)\times B(x_s,\epsilon_k)
\end{equation*}
such that 
\begin{equation*}
\norm{u_k(x_{rs}^k)-u_k(x_{sr}^k)}=\norm{x_{rs}^k-x_{sr}^k}.
\end{equation*}
Using Arzel\`a--Ascoli theorem and passing to a subsequence we may assume that maps $u_k$ converge locally uniformly to some $1$-Lipschitz map $u_0$. Observe now that 
\begin{equation*}
x_{rs}^k\text{ converges to }x_r\text{ for all }r,s=1,\dotsc,m+1.
\end{equation*}
Thus, by the locally uniform convergence, $u_0$ is an isometry on $\{x_1,\dotsc,x_{m+1}\}$. Observe that 
\begin{equation*}
\mu_{\epsilon_k}\text{ converges weakly* to }\mu_0=\sum_{i=1}^{m+1}\delta_{x_i}v_i.
\end{equation*}
Now, Lemma \ref{lem:attain} tells us that $u_0$ is an optimal potential of $\mu_0$.

Suppose that points $x_1,\dotsc,x_{m+1}$ are such that for $i\neq j$, $i,j=1,\dotsc,m$,
\begin{equation}\label{eqn:inequality}
\Big\langle\frac{x_i-x_{m+1}}{\norm{x_i-x_{m+1}}} , \frac{x_j-x_{m+1}}{\norm{x_j-x_{m+1}}}\Big\rangle< \Big\langle\frac{v_i}{\norm{v_i}},\frac{v_j}{\norm{v_j}}\Big\rangle.
\end{equation}
Then if we define $h\colon\{x_1,\dotsc,x_{m+1}\}\to\mathbb{R}^m$ by 
\begin{equation*}
h(x_{m+1})=0\text{, }h(x_i)=\norm{x_i-x_{m+1}}\frac{v_i}{\norm{v_i}}\text{ for }i=1,\dotsc,m,
\end{equation*}
then $h$ is $1$-Lipschitz. By the Kirszbraun theorem we may assume that $h$ is defined on the entire space. Moreover for 
\begin{equation*}
\pi=\sum_{i=1}^{m+1}v_i\delta_{(x_i,x_{m+1})}
\end{equation*}
we have
\begin{equation*}
\mathrm{P}_1\pi-\mathrm{P}_2\pi=\mu_0
\end{equation*}
and 
\begin{equation*}
\pi=\sum_{i=1}^{m}\frac{h(x_i)-h(x_{m+1})}{\norm{x_i-x_{m+1}}}\norm{v_i}\delta_{(x_i,x_{m+1})}
\end{equation*}
Theorem \ref{thm:ae} yields that $h$ is an optimal potential and $\pi$ is an optimal transport. It follows that 
\begin{equation*}
\norm{\mu_0}_{\mathcal{W}(\mathbb{R}^n,\mathbb{R}^m)}=\sum_{i=1}^m\norm{v_i}\norm{x_i-x_{m+1}}.
\end{equation*}
Theorem \ref{thm:ae} tells us that also
\begin{equation*}
\pi=\sum_{i=1}^{m}\frac{u_0(x_i)-u_0(x_{m+1})}{\norm{x_i-x_{m+1}}}\norm{v_i}\delta_{(x_i,x_{m+1})}
\end{equation*}
As $u_0$ is an isometry on $\{x_1,\dotsc,x_{m+1}\}$, it follows that for $i,j=1,\dotsc,m$
\begin{equation*}
\norm{h(x_i)-h(x_j)}=\norm{x_i-x_j}
\end{equation*}
which is not true, as the inequalities in (\ref{eqn:inequality}) are strict.
The obtained contradiction shows that there is no such sequence $(\epsilon_k)_{k=1}^{\infty}$, i.e. there exists $\epsilon_0>0$ such that for all $\epsilon\in (0,\epsilon_0)$ there exists a transport set with non-zero measure $\mu_{\epsilon}$ for any optimal potential of $\mu_{\epsilon}$.

By Theorem \ref{thm:balance} it follows that for such $\epsilon$ there is
is no optimal transport with absolutely continuous marginals for $\mu_{\epsilon}$.
\qed
\end{proof}

The proof of the following theorem is based on the same idea as the proof of Theorem \ref{thm:nonzero}. 
Note that we do not require below that the norms on $\mathbb{R}^n$ and on $\mathbb{R}^m$ are Euclidean. For a $1$-Lipschitz map $u\colon\mathbb{R}^n\to\mathbb{R}^m$ a leaf of $u$ is a maximal, with respect to the order induced by inclusion, set $\mathcal{S}$ such that $u|_{\mathcal{S}}$ is an isometry. A transport set is defined as a set $A$ that enjoys the property that if $x\in A$ belongs to a unique leaf of $u$, then for any $y\in\mathbb{R}^n$  such that $\norm{u(y)-u(x)}=\norm{y-x}$ there is $y\in A$.
This is to say, the leaves and transport sets are defined as in the Euclidean case.

\begin{theorem}\label{thm:affine}
Let $m\leq n$. Suppose that the norm on $\mathbb{R}^m$ is strictly convex. Suppose that $\mathcal{F}$ is a locally uniformly closed subset of $1$-Lipschitz maps of $\mathbb{R}^n$ to $\mathbb{R}^m$. Suppose that $\mathcal{F}$ has the property that for any absolutely continuous measure $\mu\in\mathcal{W}(\mathbb{R}^n,\mathbb{R}^m)$ and any $\mathcal{F}$-optimal potential $u_0$ of $\mu$ 
we have $\mu(A)=0$ for any transport set $A$ of $u_0$. Then either $m=1$ or $m>1$ and
\begin{enumerate}[i)]
\item\label{i:affineprim} $m=n$, any $u\in\mathcal{F}$ is affine, and there exists $u\in\mathcal{F}$ that is an isometry of $\mathbb{R}^n$ and of $\mathbb{R}^m$,
\item\label{i:subspace} for any absolutely continuous $\mu$, any $\mathcal{F}$-optimal potential of $\mu$ is an isometry on a maximal subspace $V\subset\mathbb{R}^n$, so that 
\begin{equation}\label{eqn:condition}
\mu(\{x\in\mathbb{R}^n\mid Px\in A\})=0\text{ for any Borel set }A\subset W;
\end{equation}
here $P$ denotes a projection onto a complement $W$ of $V$.
\end{enumerate} 
Suppose that  the norms are Euclidean. Then, if any $\mathcal{F}$-optimal potential is affine and is an isometry on a maximal subspace $V\subset\mathbb{R}^n$ such that (\ref{eqn:condition}) holds true, then $\mu(A)=0$ for any transport set of its $\mathcal{F}$-optimal potential.
\end{theorem}
\begin{proof}
Suppose that $m>1$. Choose any pairwise distinct points $x_1,x_2,x_3\in\mathbb{R}^n$ and any affinely independent $v_1,v_2,v_3\in\mathbb{R}^m$ such that $\sum_{i=1}^3v_i=0$. Let 
\begin{equation*}
\nu_0=\sum_{i=1}^3v_i\delta_{x_i}.
\end{equation*}
Then $\nu_0\in\mathcal{W}(\mathbb{R}^n,\mathbb{R}^m)$. For $\epsilon>0$ let 
\begin{equation*}
\nu_{\epsilon}=\frac1{\lambda(B(0,\epsilon)}\sum_{i=1}^3v_i\lambda|_{B(x_i,\epsilon)}
\end{equation*}
Choose respective $\mathcal{F}$-optimal potentials $u_{\epsilon}$ for $\nu_{\epsilon}$. These exist as $\mathcal{F}$ is locally uniformly closed. 
Observe that, by the assumption, $\nu_{\epsilon}(B_{\epsilon})=0$ for any Borel set $B_{\epsilon}$ consisting of trivial leaves of $u_{\epsilon}$. 
Whence, $\nu_{\epsilon}$ is concentrated on non-trivial leaves of $u_{\epsilon}$.
Let $N_{i\epsilon}$ denote the union of all non-trivial leaves that intersect $B_{i\epsilon}=\mathrm{cl}B(x_i,\epsilon)$ 
 for $i=1,2,3$ and $\epsilon>0$.
The map 
\begin{equation*}
\mathbb{R}^n\setminus B_{i\epsilon}\ni x\mapsto \sup\Big\{\frac{\norm{u(x)-u(y)}}{\norm{x-y}}\mid y\in B_{i\epsilon}\Big\}\in\mathbb{R}
\end{equation*}
is lower semi-continuous.
Note that 
\begin{equation*}
N_{i\epsilon}=\Big\{x\in\mathbb{R}^n\setminus B_{i\epsilon}\mid \sup\Big\{\frac{\norm{u(x)-u(y)}}{\norm{x-y}}\mid y\in B_{i\epsilon}\Big\}=1\Big\}\cup(B_{i\epsilon}\cap A_{\epsilon}),
\end{equation*}
where $A_{\epsilon}$ denotes the union of all non-trivial leaves of $u_{\epsilon}$ and is Borel measurable by the argument of Lemma \ref{lem:non-trivial}.
Hence $N_{i\epsilon}$ is Borel measurable. 

By the assumption,
\begin{equation*}
\nu_{\epsilon}(N_{i\epsilon})=0,
\end{equation*}
which implies, as in the proof of Theorem \ref{thm:nonzero}, that 
\begin{equation*}
\norm{u_{\epsilon}(x^{\epsilon}_{rs})-u_{\epsilon}(x^{\epsilon}_{sr})}=\norm{x^{\epsilon}_{rs}-x^{\epsilon}_{sr}}
\end{equation*}
for some points 
\begin{equation*}
(x^{\epsilon}_{rs},x^{\epsilon}_{sr})\in B(x_r,\epsilon)\times B(x_s,\epsilon), r,s=1,\dotsc,3, r\neq s. 
\end{equation*}
By the Arzel\`a--Ascoli theorem and passing to a subsequence we may assume that $u_{\epsilon}$ converges locally uniformly to some $u_0\in\mathcal{F}$, which is an $\mathcal{F}$-optimal potential of $\nu_0$ by Lemma \ref{lem:attain}.
By the uniform convergence we infer that $u_0$ is isometric on $\{x_1,x_2,x_3\}$. Let now $x_2=tx_1+(1-t)x_3$ for some $t\in (0,1)$. Then any map $f\colon\{x_1,x_2,x_3\}\to\mathbb{R}^m$ that is isometric satisfies 
\begin{equation}\label{eqn:cc}
f(tx_1+(1-t)x_3)=tf(x_1)+(1-t)f(x_3).
\end{equation}
Indeed, as $f$ is isometric,
\begin{equation*}
\norm{f(x_2)-f(x_1)}= (1-t)\norm{x_3-x_1}\text{ and }\norm{f(x_3)-f(x_2)}= t\norm{x_3-x_1}.
\end{equation*}
As $\norm{f(x_3)-f(x_1)}=\norm{x_3-x_1}$ it follows that we have equality in the triangle inequality
\begin{equation*}
\norm{f(x_3)-f(x_1)}\leq\norm{f(x_2)-f(x_1)}+\norm{f(x_3)-f(x_2)}.
\end{equation*}
By the strict convexity it follows that there is $\lambda>0 $ such that 
\begin{equation*}
f(x_2)-f(x_1)=\lambda(f(x_3)-f(x_1)). 
\end{equation*}
Taking the norms we arrive at (\ref{eqn:cc}).
A function $f$ that satisfies (\ref{eqn:cc}) may be extended to $\mathbb{R}^n$ to an affine map that has derivative of operator norm at most $m$. This follows by the Hahn--Banach theorem.
As $u_0$ is isometric on $\{x_1,x_2,x_3\}$, we infer that 
\begin{equation*}
\sum_{i=1}^3\langle u_0(x_i),v_i\rangle\ \leq \sup\Big\{\sum_{i=1}^3\langle f(x_i),v_i\rangle\mid f\colon\mathbb{R}^n\to\mathbb{R}^m \text{ is linear and }\norm{f}\leq m\Big\}
\end{equation*} 
Note now that the set of vectors $v_1,v_2,v_3$ that sum up to zero and are affinely independent is dense in the set of vectors $v'_1,v'_2,v'_3$ that sum up to zero. Moreover, $u_0$ is an $\mathcal{F}$-optimal potential for $\nu_0$. We conclude that for any $u\in\mathcal{F}$ and any vectors $v_1,v_2,v_3$ that sum up to zero there is
\begin{equation*}
\sum_{i=1}^3\langle u(x_i),v_i\rangle\ \leq \sup\Big\{\sum_{i=1}^3\langle f(x_i),v_i\rangle\mid f\colon\mathbb{R}^n\to\mathbb{R}^m\text{ is linear and }\norm{f}\leq m\Big\}
\end{equation*} 
Take now $v_2=v$, $v_1=-tv$ and $v_3=-(1-t)v$ with $t\in (0,1)$ as above and any $v\in\mathbb{R}^m$. We infer that
\begin{equation*}
\langle u(x_2)-tu(x_1)-(1-t)u(x_3) ,v\rangle\leq 0.
\end{equation*}
As this holds for any $v$ we infer that $u$ is affine. This proves part of \ref{i:affineprim}).

If $u$ is affine then there exists a subspace $V\subset\mathbb{R}^n$, possibly trivial, i.e. $V=\{0\}$, such that any set of the form 
\begin{equation}\label{eqn:leaff}
\{x\in\mathbb{R}^n\mid Px\in A\}
\end{equation} 
for a Borel measurable set $A\subset W$ is a transport set of $u$. Here $P$ denotes a projection onto a complement $W$ of $V$.  Indeed, let $V\subset\mathbb{R}^n$ be a maximal subspace such that $u|_V$ is an isometry.  Suppose that $V$ is not a leaf of $u$. Then there exists $y\notin V$ such that for all $x\in V$
\begin{equation*}
\norm{u(y)-u(x)}=\norm{y-x}.
\end{equation*}
It follows that for all non-zero $\lambda\in\mathbb{R}$
\begin{equation*}
\Big\lVert u(y)-u\Big(\frac{x}{\lambda}\Big)\Big\rVert=\Big\lVert y-\frac{x}{\lambda}\Big\rVert
\end{equation*}
for all $x\in V$. Hence for all $\lambda\in\mathbb{R}$ we have $\norm{u(\lambda y)-u(x)}=\norm{\lambda y-x}$. As $u$ is affine, it is also an isometry on $V+\mathbb{R}y$. This contradiction shows that $V$ is a leaf of $u$. Thus \ref{i:subspace}) is proven.

We shall now provide an example of a vector measure $\mu$ such that for any proper subspace $V$ and any $x_0$ there is $c>0$ such that 
\begin{equation}\label{eqn:mum}
\mu\Big(\big\{x\in\mathbb{R}^n\mid \norm{P(x-x_0)}\leq c\big\}\Big)\neq 0.
\end{equation}
Choose any $x_1,\dotsc,x_{m+1}\in\mathbb{R}^n$ affinely independent. Let $\epsilon>0$ be a number such that any set $\{y_1,\dotsc,y_{m+1}\}$, with $y_i\in B(x_i,\epsilon)$, $i=1,\dotsc,m+1$, is affinely independent. Choose vectors $v_1,\dotsc,v_{m+1}\in\mathbb{R}^m$ that add up to zero and are affinely independent. Let 
\begin{equation*}
\mu=\sum_{i=1}^{m+1}v_i\lambda|_{B(x_i,\epsilon)},
\end{equation*}
where $\lambda$ denotes the Lebesgue measure. Choose any proper affine subspace $V\subset\mathbb{R}^n$. Then $V$ intersects at most $m$ of the balls $B(x_i,\epsilon)$, $i=1,\dotsc,m+1$. So does the set 
\begin{equation*}
\big\{x\in\mathbb{R}^n\mid \norm{P(x-x_0)}\leq c\big\}
\end{equation*}
provided that $c>0$ is sufficiently small. Thus (\ref{eqn:mum}) follows. It implies, by \ref{i:subspace}), that $V=\mathbb{R}^n$. We have shown that any $\mathcal{F}$-optimal potential of $\mu$ has to be an isometry. Hence $m= n$ and the proof of \ref{i:affineprim}) is complete.

To prove the last part of the theorem, it is enough to prove that the translates of $V$ are the only leaves of an affine map. This holds true, since any point in $\mathbb{R}^n$ is covered by a translate of $V$ and the leaves of a map foliate $\mathbb{R}^n$, up to Lebesgue measure zero, if the considered norms are Euclidean, c.f. \cite{Ciosmak2}.\qed
\end{proof}

%
%

\bibliographystyle{spmpsci}      
\bibliography{refs}   
%

\end{document}